\documentclass[12pt]{article}
\usepackage{amsmath,amssymb,amsfonts,bm}
\usepackage{mathrsfs}
\usepackage[numbers,sort&compress]{natbib}

\usepackage[bookmarks=false,pdfstartview=FitH,linkbordercolor={0.5 1 1},%
citebordercolor={0.5 1 0.5},unicode,%
pagebackref,hyperindex]{hyperref}%

\catcode`@=11

\@addtoreset{equation}{section}
\catcode`@=12
\newtheorem{theorem}{Theorem}[section]
\newtheorem{definition}[theorem]{Definition}
\newtheorem{lemma}[theorem]{Lemma}
\newtheorem{stat}[theorem]{Proposition}

\newtheorem{ex}[theorem]{Example}
\newtheorem{cor}[theorem]{Corollary}
\def\qed{$\ \vbox{\hrule\hbox{\vrule height1.3ex\hskip0.8ex\vrule}\hrule}$}

\def\proof{\mbox{\it Proof.\ }}
\def\rit#1{\mbox{$\mathbb{R}^#1$}}
\def\class{\mathscr{F}}
\def\norma{\mbox{$\|\cdot\| $}}
\def\span{\mbox{$\rm span$}}
\def\qcm{\mbox{$\sigma$}}
\def\ovm{\mbox{$\chi$}}
\def\absco{\mbox{$\rm absco$}}
\def\co{\mbox{$\rm co$}}

\def\length{\mbox{$\ell$}}

\def\set1N{\mbox{$\{ 1,2,\ldots ,N\}$}}

\def\v2{\vspace{2mm}\noindent}

\begin{document}
\title
{\textbf{Estimates of Amplitudes of Transient Regimes in
Quasi--Controllable Discrete Systems}\thanks{V.Kozyakin was
partially supported by the the Russian Foundation for Fundamental
Research Grant 930100884, A. Pokrovskii was partially supported by
the Australian Research Council Grant A 8913 2609} }
\author
{
V. Kozyakin \\
{\small\sl Institute of Information}\\
{\small\sl Transmission Problems,}\\
{\small\sl Russia, Moscow, Ermolovoj st.,19.}\\
{\small\sl e-mail kozyakin@ippi.msk.su}
\and
A.Pokrovskii \thanks{{\it Permanent address:} Institute
of Information Transmission Problems, Russian Academy of Sciences. }\\
{\small\sl Mathematics Department,}\\
{\small\sl University of Queensland,}\\
{\small\sl 4072 Australia.}\\
{\small\sl e-mail ap@maths.uq.oz.au}
}
\date{}
\maketitle

\begin{abstract}
Families of regimes for discrete control systems are studied
possessing a special {\em quasi--con\-trol\-\l\-abi\-li\-ty property}
that is similar to the Kalman con\-trol\-\l\-abi\-li\-ty property.
A new approach is proposed to estimate the amplitudes of
transient regimes in quasi--controllable systems.
Its essence is in obtaining of constructive
a priori bounds for degree of overshooting in terms of the
quasi--con\-trol\-\l\-abi\-li\-ty measure.
The results are applicable for analysis of transients, classical
absolute stability problem and, especially, for
stability problem for desynchronized systems.
\end{abstract}

\v2

\v2

\section*{Introduction}
Currently, there are a growing number of cases in which systems are
described as operating permanently as if in the transient mode.
Examples are flexible manufacturing systems, adaptive control systems
with high level of external noises, so called desynchronized systems
or asynchronous discrete event systems
\cite{AKKK,KKKK,CHINA}.
In connection with this, it is necessary to ensure that
the state vector amplitude satisfies reasonable estimates within the whole
time interval of the system functioning including the interval of
transient regime and an infinite interval when the state vector is
``close to equilibrium''.
Emphasize,  that this necessity often contradicts
the usual desire to design a feedback which
makes the system as stable as possible.
The reason is that the stability property
characterizes only the asymptotic behavior of a system and does not
take into account system behavior during the transient
interval. As a result, a stable system can have  large overshooting or
``peaks'' in the transient process that can result in complete failure
of a system.
First mentions about systems with peak effects could be found in
\cite{BonYou68,BonYou70} and \cite{Mita76}. In
\cite{Kim81,Mita77,Mita78,Pol80b,Zei83} this effect was studied for
some classes of linear systems.
As it was noted in
\cite{MiYo80,Izm87a,Izm88,Izm89,Kim81,Sha84,OlCi88}
when the regulator in feedback links is chosen to guarantee as large
degree of stability as possible then, simultaneously, overshooting of
the system state during the transient process grows i.e., the peak
effects  are getting more dangerous.
From the geometrical point of view peak effect means that
when we are trying to design a feedback which improves
the stability of the system we should ``spoil'' automatically
a form of Lebesgue surfaces of respective Lyapunov functions.

The above  papers were mainly concerned with
continuous time control systems because, in completely controllable
and observable discrete system it is possible to chose
feedback which turns to zero  the specter
of the respective closed--loop system.
Nevertheless,  similar effects occurred when optimizing asymptotic
behaviour of  {\em badly controllable or observable} discrete time systems
which arise in some  applications, see further references in
\cite{Wi,Vo}.
Consider as the simplest, if trivial,  example the linear system
which is described by the relations
\begin{equation}\label{E0}
x_{n+1}=Ax_{n}+bx_{n}^{1},
\qquad
n = 0,1,2,\ldots \ .
\end{equation}
Here  $x=(x^{1},x^{2})^{T} $
be a vector from $\rit{2}$, $A$ be a matrix of the form
$$
\left(
\begin{array}{cc}
a & \varepsilon\\
\varepsilon & a
\end{array}
\right)
$$
with the small $\varepsilon$
and $b\in\rit{2}$ defines a feedback to be constructed.
From the asymptotical point of view the best vector
$b_{*}$ is $\left(-2a,-\frac{a^{2}+\varepsilon^{2}}{\varepsilon}\right)$ which
makes the eigenvalues of a closed system equal to zero.
On the other hand, for small $\varepsilon$
this vector $b_{*}$ is the most dangerous at the first time step, because
the system (\ref{E0}) can be written for this $b_{*}$ as
$x_{n+1}=A_{*}x_{n}$ where
$$
\left(
\begin{array}{cc}
-a & \varepsilon\\
\frac{a^{2}}{\varepsilon} & a
\end{array}
\right)
$$
has a big element $\frac{a^{2}}{\varepsilon}$ in the left bottom corner.
There arises a general question if this kind of the peak
effect is connected only
with poor controllability or observability of the system? If an
answer is positive, then the respective quantitative estimates
are of interest.
Especially urgent such estimates seems to be when a whole class
of systems is examined just as in problems of absolute
stability or in desynchronized systems.
Another schemes of appearing peak  effects in discrete systems
see  in \cite{DIJ,IV}.

In this paper a new approach is developed presenting the means to solve
for some classes of systems
effectively the problem of estimation the state vector amplitude within
the whole time interval. The key concept used is a
quasi--controllability property of a system that is similar to
the Kalman controllability property.
The degree of quasi--controllability can be characterized by a numeric
value. The main result of the paper is in proving the following:
if  a quasi--controllable system is stable then the amplitudes of all
its state trajectories starting from the unit ball are bounded by the
value reciprocal of the quasi--controllable measure. Due to the fact
that the measure of quasi--controllability can be easily computed, this
fact becomes an efficient tool for analysis of transients.
It is shown
also that for quasi--controllable systems the properties of stability
or instability are robust with respect to small perturbation of system's
parameters.
Some other results in this direction were announced  in \cite{KP,KKP}.

\section{Quasi--controllable families of matrices}
The notion of {\em quasi--controllability} of the system
will be introduced in this section.
Degree of quasi--controllability will be  estimated by
some nonnegative value, the {\em quasi--controllability
measure}.
The basic  property of quasi--controllability measure and
some examples will be also
discussed in this section.

\subsection{Definition and the first properties} Let ${\class}
=\{A_{1}$, $A_{2}$, \ldots , $A_{M}\}$ be a finite family of real
$N\times N$ matrices.

\begin{definition}\label{quasicon} A family ${\class}$ is said to be
quasi--controllable one if  no nonzero proper subspace of \rit{N} is
invariant for all matrices from $\class$.
\end{definition}
\noindent

Denote by ${\class}_{k}\ (k = 1,2,\ldots )$ the set of finite products
of matrices from ${\class}\bigcup\{I\}$ which contain no more that $k$
factors. Define ${\class}_{k}(x), \ x\in\rit{N}$, as the set of vectors
$Lx$, with $L\in{\class}_{k}$.  Denote by ${\co}(W)$ and ${\span}(W)$
respectively the convex and the linear hulls of the set
$W\subseteq\rit{N}$. Introduce also the set ${\absco}(W) =
{\co}(W\bigcup -W)$ which is called the {\em absolute convex hull} of
$W$.  Let {\norma} be a norm in \rit{N}; a ball in this norm of the
radius  $t$ centered at $0$ denote by ${\bf S}(t)$.

\begin{theorem}\label{L11}
Suppose that $p\ge N-1$. Then a  family ${\class}$
is quasi--controllable if and only if
${\span}\{{\class}_{p}(x)\} =\rit{N}$
for each nonzero $x\in\rit{N}$.
\end{theorem}

\proof
 Let the family ${\class}$ be quasi--controllable and
$x\in\rit{N}$ be a given nonzero vector. Introduce the sets
$\mathscr{L}_{0} = {\span}\{x\}$ and $\mathscr{L}_{k} =
{\span}\{{\class}_{k}(x)\}$, $k\ge 1$. Then
\begin{equation}\label{E11}
\mathscr{L}_{0}\subseteq\mathscr{L}_{1}\subseteq \ldots \subseteq
\mathscr{L}_{p}\subseteq\rit{N}.
\end{equation}
Therefore,
\begin{equation}\label{E12}
1\le\dim\mathscr{L}_{0}\le\dim\mathscr{L}_{1}\le\ldots\le\dim\mathscr{L}_{p}\le
N.
\end{equation}
On the other hand,
\begin{equation}\label{E13}
A_{i}\mathscr{L}_{j}\subseteq\mathscr{L}_{j+1},\qquad A_{i}\in
{\class},0\le j\le p-1.
\end{equation}

If $\dim\mathscr{L}_{p} = N$ then $\mathscr{L}_{p} =
{\span}\{{\class}_{p}(x)\} =\rit{N}$. If $\dim\mathscr{L}_{p}<N$
then by (\ref{E12}) and the condition
 ${p\ge N-1}$, the equality $\dim\mathscr{L}_{j}=\dim\mathscr{L}_{j+1}$ holds
for some ${j\in [0,p-1]}$.  The last equality and (\ref{E11}) imply
$\mathscr{L}_{j} = \mathscr{L}_{j+1}$.  By the last equality and
(\ref{E13}) the subspace $\mathscr{L}_{j}$ should be invariant with
respect to all matrices from $\class$; due to quasi--controllability
of the family ${\class}$ this subspace coincides with \rit{N}.
Hence, $\mathscr{L}_{j} = \mathscr{L}_{j+1} = \ldots =
\mathscr{L}_{p} = {\span}\{{\class}_{p}(x)\} =\rit{N}$.

\vspace{2mm} \noindent Now suppose that  ${\span}\{{\class}_{p}(x)\}
=\rit{N}$, but the family ${\class}$ is not quasi--controllable.
Then there exists a nonzero proper subspace $\mathscr{L}\subset
\rit{N}$ which is invariant with respect to all matrices from
$\class$.  In this case the inclusion
${\span}\{{\class}_{p}(x)\}\subseteq \mathscr{L}$ holds for each
${x\in \mathscr{L}}$.  Therefore,
${\span}\{{\class}_{p}(x)\}\neq\rit{N}$. This contradiction proves
the quasi--controllability of the family ${\ class}$. The lemma is
proved. \qed

\begin{definition}\label{pmesquasi}
The value  {\qcm}$_{p}({\class})$ defined by
$$
{\qcm}_{p}({\class}) = \inf_{x\in\rit{N},\|x\|=1}
\sup\{t:{\bf S}(t)\subseteq{{\absco}[\class}_{p}(x)]\}
$$
is called $p$-measure of quasi--controllability of the family ${\class}$
(with respect to the norm {\norma}).
\end{definition}

\begin{theorem}\label{L12} Suppose that  ${p\ge N-1}$. The family ${\class}$
is quasi--controllable if and only if ${{\qcm}_{p}({\class})\neq 0}$.
\end{theorem}

\proof
Suppose that ${{\qcm}_{p}({\class})\neq 0}$. Then
${\bf S}[\|x\|{\qcm}_{p}({\class})]\subseteq{{\absco}[\class}_{p}(x)]$
holds
for each nonzero ${x\in\rit{N}}$
and, further,
$\rit{N}={\span}\{{\class}_{p}(x)\}$. Therefore, by Theorem
\ref{L11} the family ${\class}$ is quasi--controllable.

\vspace{2mm}
Suppose now that the family
${\class}$ is quasi--controllable but
${\qcm}_{p}({\class})=0$. Then there exist $x_{n}\in\rit{N}$, ${\|x_{n}\|= 1}$,
and $y_{n}\in{{\absco}[\class}_{p}(x_{n})]$ such that $y_{n}\rightarrow 0$
and
$ty_{n}\not\in{{\absco}[\class}_{p}(x_{n})]$ for $t>1$.  Without loss
of generality we can suppose that the sequences
$\{x_{n}\}$ and
$\{\frac{y_{n}}{\|y_{n}\|}\}$ are convergent: $x_{n}\rightarrow x$,
$\frac{y_{n}}{\|y_{n}\|}\rightarrow z$.

By Theorem \ref{L11} the linear hull of the set $\{{\class}_{p}(x)\}$
coincides with \rit{N}. Hence, there exist matrices $L_{1}$, $L_{2}$, \ldots ,
$L_{N}\in{\class}_{p}$ such that the vectors $L_{1}x$, $L_{2}x$, \ldots ,
$L_{N}x$ are linearly independent. Then
the vectors
$L_{1}x_{n}$, $L_{2}x_{n}$, \ldots , $L_{N}x_{n}$
are also independent for all sufficiently large $n$.
It means that for any $n$ there exist numbers
$$
\theta ^{(n)}_{1}, \theta ^{(n)}_{2}, \ldots , \theta ^{(n)}_{N},\qquad
\sum^{N}_{i=1}\theta ^{(n)}_{i}=1\ ,
$$
such that the vector
\begin{equation}\label{E14}
z_{n} =\sum^{N}_{i=1}\theta ^{(n)}_{i}L_{i}x_{n}
\end{equation}
is collinear to  $y_{n}$ i.e.,
$z_{n}=\eta _{n}y_{n}$ $(\eta _{n}>0)$.

By definition
$z_{n}\in{\absco}\{L_{1}x_{n}$, $L_{2}x_{n}$, \ldots ,
$L_{N}x_{n}\}\subseteq{\absco}[{\class}_{p}(x_{n})]$ and $ty_{n}$
does not belong to the set
${\class}_{p}(x_{n})$ for  $t>1$. Therefore,
$\eta _{n}\le 1$. The last inequality and the condition $y_{n}\rightarrow 0$
imply
$z_{n}\rightarrow 0$.
Without loss of generality the sequences
$\{\theta ^{(n)}_{1}\}$, $\{\theta ^{(n)}_{2}\}$, \ldots ,
$\{\theta ^{(n)}_{N}\}$ can be supposed to be convergent
to some limits
$\theta _{1}$, $\theta _{2}$, \ldots , $\theta _{N}$.
Now, after transition to the limit in
(\ref{E14}), we get
$$
\sum^{N}_{i=1}\theta _{i}L_{i}x=0,\qquad \sum^{N}_{i=1}\theta _{i}=1\ .
$$
This contradicts the linear independence of the vectors $L_{1}x$,
$L_{2}x$, \ldots , $L_{N}x$, and the theorem is proved. \qed

The following theorem is useful
when a family of matrices depends on a parameter.
\begin{theorem}\label{L13}
Let ${p\ge N-1}$ and  the $N\times N$ matrices
$$
A_{1}(\tau),A_{2}(\tau ), \ldots , A_{M}(\tau)
$$
be continuous at
the point  $0$ with respect to the real parameter
$\tau$. Suppose that the family
${\class}(\tau) =\{A_{1}(\tau)$, $A_{2}(\tau)$, \ldots , $A_{M}(\tau)\}$
is quasi--controllable at $\tau = 0$.
Then the family
${\class}(\tau )$ is quasi--controllable
for all sufficiently small $\tau$
and the function
${\qcm}_{p}[{\class}(\tau)]$ is continuous in $\tau$
at the point  ${\tau=0}$.
\end{theorem}

The proof is relegated to the Appendix.

\subsection{Examples}
Let $A$ be a matrix of the size $N$ and $b,c\in\rit{N}$.
Consider the family ${\class}={\class}(A,b,c)$ which consists
of the matrix $A$ and the matrix $Q=b\,c^{T}$ with elements
$q_{ij}=b_{i}c_{j}$, $i,j=1,\ldots,N$.
\begin{stat}\label{Ex10}
The family ${\class}(A,b,c)$ is quasi--controllable
if and only if the pair $(A,b)$ is completely controllable
and the pair $(A,c)$ is completely observable.
\end{stat}

\proof
Evidently, the subspace $E\subset \rit{N}$ is invariant with
respect to the matrix $Q$ if and only if either $b\in E$ or
$E\subset c^{0}$ where
$$
c^{0}=\left\{x\in\rit{N}:\: \sum_{i=1}^{N}x^{i}c^{i}=0\right\}.
$$
Farther,
the matrix $A$ has a proper invariant subspace $E_{1}$
which contains the vector $b$ if and only if
$$
\span(\{b,Ab,\ldots,A^{N-1}b\})=\rit{N},
$$
that is if the pair $(A,b)$ is completely controllable.
At last,
the matrix $A$ has a proper invariant subspace $E_{2}$
which is contained in $c^{0}$ if and only if
$$
\span(\{c,cA,\ldots,cA^{N-1}\})=\rit{N},
$$
that is if the pair $(A,c)$ is completely observable.
Therefore the assertion is proved. \qed

\vspace{2mm}
\noindent
The following example is the most important for this paper.
Let us consider a $N\times N$ scalar matrix $A = (a_{ij})$
of the size $N$ and introduce the family
${\class}_{1}(A) =\{A_{1}$, $A_{2}$, \ldots , $A_{N}\}$
by equalities
\begin{equation}\label{mixture}
A_{i} =\left(\begin{array}{cccccc}
1 & 0 &\dots & 0 &\dots & 0\\
0 & 1 &\dots & 0 &\dots & 0\\
\vdots &\vdots &\ddots &\vdots &\ddots &\vdots\\
a_{i1}& a_{i2}&\dots & a_{ii}&\dots & a_{iN}\\
\vdots &\vdots &\ddots &\vdots &\ddots &\vdots\\
0 & 0 &\dots & 0 &\dots & 1\end{array}\right).
\end{equation}
The families ${\class}_{1}(A)$ play a key role in the theory
of desynchronized systems, see \cite{AKKK,KKKK}
and, also, Subsection \ref{ex2S}.

The matrix $A$ is said to be {\em irreducible}, if by any reordering of the
basis elements in $\rit{N}$ it cannot be represented in a block triangle form.
$$
A=\left(\begin{array}{cc} B&C\\0&D\end{array}\right).
$$
Irreducibility  of the matrix $A$ means that this matrix
has no nonempty proper invariant subspace
which is the linear hull of a subset of the basic vectors
$$
e_{i} =(0,0, \ldots,1, \ldots , 0)\qquad i = 1,2,\ldots ,N.
$$

Let the norm ${\norma}$ ¢ $\rit{N}$ is defined by
$\|x\| =|x_{1}| + |x_{2}| +\ldots + |x_{N}|$. Let
$$
\alpha ={\frac{1}{2N}}\min\{\|(A-I)x\|:\|x\|=1\}
\qquad
\beta ={\frac{1}{2}}\min\{|a_{ij}| : i\neq j, a_{ij}\neq 0\},\qquad
$$
\begin{stat}\label{Ex11} The family
${\class}_{1}(A)$ is quasi--controllable,
if and only if $1$ is not an eigenvalue of
$A$ and the matrix
$A$ is irreducible. If
${\class}_{1}(A)$ is quasi--controllable then
$$
{\qcm}_{N}[{\class}_{1}(A)]\ge {\alpha}{\beta}^{N-1}.
$$
\end{stat}

\proof Let 1 be an eigenvalue of $A$ with an eigenvector
$x_{*}$. Then $x_{*}$ is an eigenvector with the eigenvalue 1
for each matrix
$A_{1}$, $A_{2}$, \ldots , $A_{N}$.
Hence, in this case the family ${\class}_{1}(A)$ is not
quasi--controllable.

Suppose that the matrix $A$ is irreducible.
Then we can assume without loss of generality
that some subspace of the form
$E_{p} = {\span}\{e_{1}$, $e_{2}$,
\ldots , $e_{p}\}$ with $p<N$ is invariant with respect
to the matrix $A$. Therefore,
$E_{p}$ should be also invariant with respect
to each matrix
$A_{1}$, $A_{2}$, \ldots , $A_{N}$.
That is, the family ${\class}_{1}(A)$ is not quasi--controllable.

Let us  now prove that the family ${\class} = {\class}_{1}(A)$
is quasi--controllable, providing that 1 is not  an eigenvalue
of  $A$ and  that  $A$ is irreducible. It will suffice to show that,
for each nonzero vector $x\in\rit{N}$,
\begin{equation}\label{E119}
{\span}\{{\class}_{N}(x)\} =\rit{N}.
\end{equation}
Choose a vector $x\in\rit{N},\|x\| = 1$, and consider the vectors
$(A_{1}-I)x$, $(A_{2}-I)x$, \ldots ,$(A_{N}-I)x\in{{\span}\{{\class}_{1}(x)\}}$.
By definition
$$
(A-I)x = (A_{1}-I)x + (A_{2}-I)x +\ldots + (A_{N}-I)x \ ,
$$
and 1 is not an eigenvector of the matrix A.
Therefore, at least one of the vectors
$(A_{1}-I)x$, $(A_{2}-I)x$, \ldots , $(A_{N}-I)x$ is nonzero.
Without loss of generality we can assume that
$(A_{1}-I)x\neq 0$ and
$\|(A_{1}-I)x\|\ge \frac{1}{N}\|(A-I)x\|\ge 2{\alpha}$.
But
\begin{equation}\label{E120}
(A_{i}-I)x = {\langle\tilde{a}_{i},x\rangle}e_{i},
\qquad i = 1,2,\ldots ,N,
\end{equation}
where ${\langle \cdot ,\cdot\rangle}$ denotes the inner
product in
\rit{N} and the vectors  $\tilde{a}_{i}$ are of the form
$$
\tilde{a}_{i} =
(a_{i1},a_{i2}, \ldots , a_{ii}-1, \ldots , a_{iN}), \qquad
i = 1,2,\ldots ,N \ .
$$
Hence,
${\langle\tilde{a}_{1},x\rangle}e_{1}\neq 0$,
${\langle\tilde{a}_{1},x\rangle}e_{1}\in{\span}\{{\class}_{1}(x)\}$
and
$\|{\langle\tilde{a}_{1},x\rangle}e_{1}\|\ge 2{\alpha}$.
This implies that
\begin{equation}\label{E121}
e_{1}\in{\span}\{{\class}_{1}(x)\}
\end{equation}
So the vector
$$
\frac{1}{2}{\langle\tilde{a}_{1},x\rangle}e_{1} =
\frac{1}{2}A_{1}x-\frac{1}{2}x
$$
belongs to ${\absco}\{{\class}_{1}(x)\}$, and, further,
${\absco}\{{\class}_{N}(x)\}$.
Therefore,
$$
{\alpha}e_{1}\in{\absco}\{{\class}_{N}(x)\}.
$$
Let $Ae_{1}= (v_1, v_2, \ldots ,v_N )$. By irreducibility of the matrix
$A$, the subspace ${\span}\{e_{1}\}$ is noninvariant with respect to
$A$. So, at least one of the coordinates $v_1$, $v_2$, \ldots ,$v_N$ of
the vector $Ae_{1}$, with the index different from 1, is nonzero.
Without loss of generality, assume that $v_{2}\not = 0$.  But the
second coordinate of the vector $Ae_{1}$ coincides with the second
coordinate of the vector $A_{2}e_{1}$ and, consequently, of the vector
$(A_{2}-I)e_{1}$.  That is, $(A_{2}-I)e_{1}\neq 0$ and, by
(\ref{E121}), $(A_{2}-I)e_{1}\in {\span}\{{\class}_{2}(x)\}$.
Therefore, by (\ref{E120})
$$
e_{2}\in{\span}\{{\class}_{2}(x)\}
$$
and the vector
$\frac{1}{2}a_{21}e_{2} =
\frac{1}{2}{\langle\tilde{a}_{2},e_{1}\rangle}e_{2} =
\frac{1}{2}A_{2}e_{1}-\frac{1}{2}e_{1}$
belongs to
${\absco}\{{\class}_{2}(e_{1})\}$
and, further, belongs to
${\absco}\{{\class}_{N}(e_{1})\}$. Hence,
$$
{\alpha}{\beta}e_{2}\in{\absco}\{{\class}_{N}(x)\}.
$$
Similarly, the irreducibility of the matrix
$A$ implies the inclusions
\begin{equation}\label{E123}
e_{i}\in{\span}\{{\class}_{i}(x)\},\quad
{\alpha}{\beta}^{i-1}e_{i}\in{\absco}\{{\class}_{N}(x)\},
\qquad i=1,2,\ldots ,N
\end{equation}
for an appropriate
reordering of the basic vectors
$e_{1}$, $e_{2}$, $e_{3}$, \ldots , $e_{N}$.
The equality (\ref{E119}) and the estimate
$$
{\qcm}_{N}[\mathscr{P}_{1}(A)]\ge{\alpha}{\beta}^{N-1}
$$
follow from the relations
(\ref{E123}).
The proof of the assertion is completed.
\qed

\section{Quasi--controllability and the peak effect}
This section contains the main results of the paper.
We investigate the influence of
quasi--controllability on stability, instability
and transient processes of dynamical systems generated by nonautonomous linear difference equations
\begin{equation}\label{E124}
x(n+1) = A(n)x(n).
\end{equation}
A conceptually simple and effective  method to
estimate norms of solutions of difference equations
uniformly for all $n=0,1,2,\ldots$
will be described.

\subsection{A priori estimate of oversooting measure}
\begin{definition}\label{absstab}
Let  $\class$ be a family of $N\times N$ matrices. The difference
equation {\rm (\ref{E124})} is {\em Lyanunov absolutely stable} with
respect to the family ${\class}$, if there exists $\mu<\infty $,
such that for each  sequence $A(n)\in {\class}$ any solution $x(n)$
of the corresponding equation satisfies the estimate
\begin{equation}\label{E125}
\sup_{n\ge 0}\|x(n)\|\le\mu\|x(0)\|.
\end{equation}
\end{definition}


\begin{definition}\label{overmes} The smallest $\mu $ for which the estimate
{\rm (\ref{E125})} holds is called the overshooting measure
of the equation {\rm (\ref{E124})} with respect to the family
${\class}$, and is denoted by ${\ovm}({\class})$.
\end{definition}

From definition it follows that ${\ovm}({\class})$ coincides with the
smallest $\mu$, for which the estimate (\ref{E125}) holds with
respect to all solutions (\ref{En124}).

\begin{theorem}\label{T11} Let the equation {\rm (\ref{E124})} be Lyapunov
absolutely stable
with respect to the  quasi\---controll\-able family
${\class}$. Then the inequality
\begin{equation}\label{mainE}
{\ovm}({\class})\le{\qcm}_{p}^{-1}({\class})
\end{equation}
holds for each $p\ge N-1$.
\end{theorem}

This assertion is the central result of the paper.
The  proof is relegated to the next subsection. Now
we will discuss some applications of the inequality (\ref{mainE}).
Clearly, the Lyapunov absolute stability of the equation (\ref{E124})
is equivalent to the Lyapunov stability of the difference
inclusion
\begin{equation}\label{En124}
x(n+1) \in F_{\class}x(n).
\end{equation}
where $F_{\class}$ is defined by
$$
F_{\class}(x)=\co\{Ax:A\in\class\}.
$$
Inclusions of the form (\ref{En124}) embrace the usual systems
of the discrete absolute stability theory \cite{KrPo,Li,NC}.
On the other hand, the Lyapunov absolute stability follows from the
absolute stability of the corresponding system. Consequently,
when estimating  overshooting measure
of control systems, it is possible to combine the classical methods of absolute stability
theory with Theorem \ref{T11}.
A definitive example of using this approach
will be presented in Subsection \ref{ex2S}.
Now let us give only some simple corollaries of Theorem
\ref{T11}.

Consider a difference equation
\begin{equation}\label{naE}
x_{n+1}=Ax_{n}+bu_{n},
\qquad
n=0,1,\ldots ,
\end{equation}
with $b\in\rit{N}$ and the scalars $u_{n}$ satisfying
for a fixed $c\in\rit{N}$ the inequality
$$
|u_{n}|\le\gamma\langle c,x_{n}\rangle
$$
where $\gamma$ is a real parameter. Such equations are common in control theory \cite{Li}

\begin{cor}\label{1aC}
Let the pair $(A,b)$ be completely controllable and the pair
$(A,c)$ be completely observable and suppose that
$\max_{|\omega|=1}(c,(\omega I-A)^{-1}b)<1$.
Then for each $p\ge 1$ any solution  $x_{n}, n=0,1,\ldots$
of the equation {\rm (\ref{naE})}
satisfies the inequality
$\|x_{n}\|\le{\qcm}_{p}^{-1}({\class}_{*})\|x_{0}\|$
where ${\class}_{*}=\{A-\gamma\,bc^{T}, A+\gamma\,bc^{T}\}$.
\end{cor}

\proof
By virtue of the proposition \ref{Ex10} the class ${\class}_{*}$
is quasi--controllable.
So this corollary follows immediately from Theorem \ref{T11} and
the circle criteria of absolute stability, \cite{Li}. \qed

Consider again the general inclusion (\ref{En124}).

\begin{cor}\label{2aC}
Let the family $\class$ be quasi--controllable and suppose that
each uniformly bounded solution
$\ldots,x_{-n},\ldots, x_{-2},x_{-1},x_{0}$
is the zero solution.
Then for each $p\ge 1$ any solution  $x_{n}, n=0,1,\ldots$
of the inclusion {\rm (\ref{En124})}
satisfies the inequality
$\|x_{n}\|\le{\qcm}_{p}({\class})\|x_{0}\|$.
\end{cor}

\proof
This corollary follows  from Theorem \ref{T11} and
from the principle of absence of any bounded solution
in the absolute stability problem \cite{KrPo}. \qed

\subsection{Proof of Theorem \protect\ref{T11}}\label{prS}
Firstly, let us
establish two auxiliary assertions. Let  $\mathscr{R}$ denote the
set of all finite products of matrices from $\class$. Define the
{\em length} $\length(R)$ of a matrix
 $R\in\mathscr{R}$
as the smallest number of factors $A_{1},A_{2},\ldots ,A_{q}\in{\class}$
in the representation $R = A_{1}A_{2}\ldots A_{q}$.

\begin{lemma}\label{increase}
Let the family ${\class}$ be quasi--controllable
and suppose that the inequalities
\begin{equation}\label{EL1}
\|Rx_*\| > \mu\frac{1}{{\qcm}_{p}({\class})}\|x_*\|,\qquad \mu >1
\end{equation}
hold for some
 $x_*\in\rit{N}$ $(x_* \not = 0)$, $p\ge N-1$, $R\in\mathscr{R}$.
Then for any $x\in\rit{N}$, $x\not = 0$ there exists a matrix
$R_x\in\mathscr{R}$ such that $\|R_{x}x\|\ge\mu\|x\|$, and
$\length(R_{x})\le \length(R)+p$.
\end{lemma}

\proof
Let us fix an arbitrary $x\in\rit{N}$, $x\not = 0$.
The vector ${\qcm}_{p}({\class})x_*$
belongs to the absolute convex hull of the set
${\class}_{p}(\frac{\|x_*\|}{\|x\|}x)$
by the definition of the quasi--controllability measure.
Therefore, there exist scalars
 $\theta_{1}$, $\theta _{2}$, \ldots , $\theta _{Q}$ with
\begin{equation}\label{E127}
\sum^{Q}_{i=1}|\theta _{i}|\le 1,
\end{equation}
and  matrices
$L_{1},L_{2},\ldots ,L_{Q}\in{\class}_{p}$
such that
$$
\sum^{Q}_{i=1}\theta _{i}\frac{\|x_{*}\|}{\|x\|}L_{i}x
={\qcm}_{p}({\class})x_{*}.
$$
Hence,
$$
\sum^{Q}_{i=1}\theta _{i}RL_{i}x
={\qcm}_{p}({\class})\frac{\|x\|}{\|x_{*}\|}Rx_{*},
$$
and, further, by  (\ref{EL1}),
$$
\sum^{Q}_{i=1}\|\theta _{i}L_{i}Rx\|\ge\mu\|x\|.
$$
But then  (see (\ref{E127})) there exists an index $i$, $1\le i\le Q$
such that the matrix
$R_{x}=L_{i}R$ satisfies $\|R_{x}x\|\ge\mu\|x\|$.

It remains to note that
the length
$\length(R_{x})\le \length(R)+p$,
due to the inclusion $L_{i}\in{\class}_{p}$,
and the lemma is proved. \qed

\begin{definition}\label{expunst}
The equation {\rm (\ref{E124})} is said to be
absolutely exponentially unstable
with degree $\lambda >1$ in the family ${\class}$,
if for some $\kappa >0$ and for each vector $x\in\rit{N}$, $x\not = 0$,
there exists a sequence $A(n)\in\class$, such that the solution
$x(n)$ of the equation {\rm (\ref{E124})} with the initial condition
$x(0)=x$ satisfies the estimate
\begin{equation}\label{Eexp}
\|x(n)\|\ge\kappa\lambda^{n}\|x(0)\|,\qquad n=0,1,2,\ldots .
\end{equation}
\end{definition}

\begin{lemma}\label{expo}
Let the family ${\class}$ be bounded and
suppose that the conditions of Lemma \ref{increase} hold. Then
the equation {\rm (\ref{E124})} is absolutely exponentially unstable in
the family ${\class}$.
\end{lemma}

\proof
Let us fix an arbitrary vector $x\in\rit{N}$, $x\not = 0$,
and construct an auxiliary sequence of vectors
$\{z(m)\}$,
$m=0,1,\ldots $, by relations  $z(0)=x$ and
$$
z(m)=R_{z(m-1)}z(m-1),\qquad m=1,2,\ldots \ .
$$
Here $R_{z(m)}$ are the matrices from Lemma \ref{increase}.
Then by Lemma \ref{increase}
\begin{equation}\label{C1}
\|z(m)\|\ge\mu^m\|z(0)\|,\qquad m=0,1,2,\ldots \ .
\end{equation}

By definition, matrices $R_{z(m)}$, $m=0,1,\ldots $, can be
represented in the form
$$
R_{z(m)}=A_{m,l(m)},\ldots ,A_{m,2},A_{m,1},\qquad A_{m,j}\in\class,
$$
where $l(m)$ is the length of $R_{z(m)}$. Denote by $\{A(n)\}$,
$n=0,1,\ldots $, the sequence of matrices
$$
A_{0,1},A_{0,2},\ldots ,A_{0,l(0)},A_{1,1},A_{1,2},\ldots ,A_{1,l(1)},\ldots
,A_{m,1},A_{m,2},\ldots ,A_{m,l(m)},\ldots \ ,
$$
and consider the solution $x(n)$ of the respective equation (\ref{E124}),
with the initial condition $x(0)=x$. Then the relations
$$
x(q_{m})=z(m),\qquad m=0,1,\ldots ,
$$
hold with $q_{0}=0$ and
$$
q_{m}=\sum^{m-1}_{i=0}l(i),\qquad m=1,2,\ldots .
$$
Estimates (\ref{C1}) imply
\begin{equation}\label{C2}
\|x(n)\|\ge\mu^{m}\|x(0)\|,\quad n=q_{m},\qquad m=0,1,\ldots \ .
\end{equation}

Norms of matrices from $\class$ are uniformly bounded by the
conditions of the lemma and also the estimates
\begin{equation}\label{C3}
q_{m}-q_{m-1}=l(m-1)\le K,\qquad m=1,2,\ldots ,
\end{equation}
hold by Lemma \ref{increase}.
Therefore, the inequality (\ref{C2}), in a  slightly weaker form, can
be extended on the positive integers
$n$ from the interval $(q_{m-1},q_{m}]$:
\begin{equation}\label{C4}
\|x(n)\|\ge\nu\mu^{m}\|x(0)\|,\qquad \nu >0,\ q_{m-1}<n\le q_{m},\
m=0,1,\ldots \ .
\end{equation}
Inequalities (\ref{C4}) for appropriate
$\kappa>0$, $\lambda >1$ imply the estimate
(\ref{Eexp}), taking into account that
$q_{m}\le mK$, $m=0,1,\ldots $,
by virtue of  (\ref{C3}).
Therefore, the lemma is proved. \qed

Let us return to and finish the proof of Theorem \ref{T11}.
Suppose that the theorem is false. Then there exists a sequence of
matrices
$\{A(n)\in\class ,\ n=0,1,\ldots\ \}$ and a solution $x(n)$
of the corresponding equation (\ref{E124}), such that
\begin{equation}\label{auE}
\|x(n_{0})\| > {\qcm}_{p}^{-1}({\class})\|x(0)\|.
\end{equation}
holds for some
$n_{0}\ge 1$, $p\ge{N-1}$.
The inequality (\ref{auE}) implies
$$
\|A(n_{0}-1)\ldots A(1)A(0)x(0)\| > {\qcm}_{p}^{-1}({\class})\|x(0)\|.
$$
Hence, by Lemma \ref{expo}, the equation (\ref{E124}) is absolutely
exponentially unstable wit respect to the family
${\class}$ and yet this equation is not even
Lyapunov absolutely stable with respect to this family.
This contradiction proves the theorem. \qed

\subsection{Application to desynchronized systems}\label{ex2S}
Recently much attention was paid to the development of methods for the
analysis of dynamics of multicomponent systems with asynchronously
interacting subsystems (see \cite{KKKK,CHINA} for further
references).  As examples we can mention the systems with faults in data
transmission channels, multiprocessor computing and telecommunication
systems, flexible manufacturing systems and so on. It turned out that under
weak and natural assumptions systems of this kind possess strong
properties like robustness. In applications the robustness is often treated as
reliability of a system with respect to perturbations of various kinds such as drift of parameters, malfunctions or noises in data transmission channels,
etc.

Let us introduce basic notions of the desynchronized systems theory.
Consider a linear system $S$ consisting of $N$ subsystems $S_{1},$ $S_{2},$
\ldots ,$S_{N}$ that interact at some discrete instants $\{T^n\}$, $-\infty
<n<\infty $.  The interaction times may be chosen according to some
deterministic or stochastic law but generally they are not known in advance.
Let the state of each subsystem $S_{i}$ be determined within the interval
${[T^{n},T^{n+1})}$ by a numerical value $x_{i}(n)$, $-\infty <n<\infty $.

Suppose that at each instant $T^{n}\in\{T^{k}:\ \ -\infty <k<\infty \}$ only one of the
subsystems $S_{i}$, $i=i(n)\in\{1,2,\ldots ,N\}$, may
change its state and the law of the state updating is linear:
$$
x_{i}(n+1)=\sum_{i=1}^{N}a_{ij}x_{j}(n),\qquad i=i(n).
$$
Consider the matrix $A=(a_{ij})$ and introduce for each $i=1,2,\ldots ,N$ an
auxiliary matrix $A_{i}$ ($i${\em -mixture} of the matrix $A$) that is
obtained from $A$ by replacing its rows with indexes $i\not =j$ with the
corresponding rows of the identity matrix $I$ (see (\ref{mixture})).  Then the dynamics equation
for the system $S$ can be written in the following compact form:

\begin{equation}\label{unrel}
x(n+1)=A_{i(n)}x(n),\qquad -\infty <n<\infty .
\end{equation}
The system described above is referred to as {\em the linear desynchronized
or asynchronous system}.

\begin{theorem}\label{Tdes}
Suppose that
$1$ is not an eigenvalue of
$A$ and that the matrix
$A$ is irreducible.
Suppose that the desynchronized system  is
Lyapunov absolutely stable.
Then
\begin{equation}\label{desE}
{\ovm}({\class})\le\frac{1}
{{\alpha}{\beta}^{N-1}}
\end{equation}
\end{theorem}

\section{Robustness of instability}
In this short concluding section we will consider another application
of the above methods to qualitative analysis of discrete systems.

Consider the difference equation
(\ref{E124}), where matrices $A(n)$ belong to a family
$\class (\tau )=\{A_{1}(\tau ),A_{2}(\tau ),\ldots ,A_{M}(\tau )\}$,
which depends on a real parameter $\tau$.

\begin{theorem}\label{unst}
Let the family $\class (\tau )$ be quasi--controllable and be
continuous at $\tau = 0$. Suppose that the equation {\rm (\ref{E124})}
is not Lyapunov absolutely stable with respect to the family $\class (0)$.
Then the equation {\rm (\ref{E124})} is not Lyapunov absolutely
stable, and, in fact, is absolutely exponentially unstable, with
respect to the family $\class (\tau )$,
for all sufficiently small $\tau$.
\end{theorem}

\proof
Suppose that the equation (\ref{E124}) is not Lyapunov absolutely
stable with respect to the family
$\class (0)$. Then there exist matrices $A(n,\tau ) \in
\class (\tau )$, $n = 0,1,\ldots $, such that at $\tau=0$
the solution of the respective  equation
(\ref{E124})
satisfies for some
 $n_{0} > 0$ the inequality
$$
\|x(n_{0})\| > {\qcm}_{N-1}^{-1}[{\class}(0)]\|x(0)\|.
$$
Therefore
\begin{equation}\label{au2E}
\|A(n_{0}-1,0)\ldots A(1,0)A(0,0)x(0)\| >
{\qcm}_{N-1}^{-1}[{\class}(0)]\|x(0)\|.
\end{equation}
On the other hand, the matrices
$\{A(n,\tau )\}$ and, by Theorem \ref{L13}, the functions
${\qcm}_{N-1}[{\class}(\tau )]$
are continuous at the point
$\tau =0$.
Consequently, (\ref{au2E}) implies
$$
\|A(n_{0}-1,\tau )\ldots A(1,\tau )A(0,\tau )x(0)\| >
{\qcm}_{N-1}^{-1}[{\class}(\tau )]\|x(0)\|.
$$
and, by virtue of Lemma \ref{expo},
the equation
(\ref{E124}) is absolutely exponentially unstable
with respect to the class
${\class}(\tau )$.  Hence, the theorem is proved. \qed

In some situations the following
corollary from the theorems \ref{L13},
\ref{T11} and \ref{unst} is useful.

\begin{cor}\label{limstab}
Let a quasi--controllable
family of matrices $\class$ =
$\{A_{1}$, $A_{2}$,\ldots , $A_{M}\}$ be the limit
of families
$\class_{m}$ = $\{A_{1,m}$, $A_{2,m}$,\ldots , $A_{M,m}\}$. Suppose
that the equation
{\rm (\ref{E124})} is Lyapunov absolutely stable with respect to families
$\class_{m}$,
$m=1,2,\ldots$. Then this equation is Lyapunov absolutely stable
with respect to the family $\class$. More than that, the families
$\class_{m}$ are quasi--controllable and the  measures
of overshooting
${\ovm}(\class_{m})$ are uniformly bounded.
\end{cor}

The following two examples show that
the previous corollary turns out to be false without the assumption
about quasi--controllability of the family
${\class}$.

\begin{ex}\label{limexp}
Consider the sequence of families $\mathscr{E}_{m}$ = $\{E_{m}\}$,
each of which consists of the single matrix
$$
E_{m}=\left(\begin{array}{cc}
1-\frac{1}{m}&1\\0&1-\frac{1}{m}
\end{array}\right).
$$
Then the limit family $\mathscr{E}$ consists of the matrix
$$
E=\left(\begin{array}{cc} 1&1\\0&1\end{array}\right),
$$
and is not quasi--controllable. Therefore the respective equation
{\rm (\ref{E124})} is not exponentially stable with respect to the
family $\mathscr{E}$, notwithstanding this equation is exponentially
stable with respect to the families $\mathscr{E}_{m}$.
\end{ex}

\begin{ex}\label{ovmunbou} Consider the sequence of the families
$\mathscr{F}_{m}$ = $\{F_{m}\}$, each of which consists of  the
single matrix
$$
F_{m} =\left(\begin{array}{cc}
1-\frac{1}{m^2}&\frac{1}{m}\\
0&1-\frac{1}{m^2}
\end{array}\right).
$$
The respective limit family $\mathscr{F}$ includes only the identity
matrix $I$ and, therefore, is not quasi--controllable. Evidently,
the equation {\rm (\ref{E124})} is stable  with respect to the
family $\mathscr{F}$, as well as with respect to the famalies
$\mathscr{F}_{m}$, $m=1,2,\ldots$. On the other hand, the  measures
of overshooting ${\ovm}(\mathscr{F}_{m})$ are not uniformly bounded.
\end{ex}

 \section*{Appendix. Proof of Theorem {\protect{\ref{L13}}}}

\proof
Establish first that, under conditions of the theorem,
there exists
${\kappa > 0}$ satisfying
\begin{equation}\label{E15}
{\qcm}_{p}[{\class}(\tau)]\ge\kappa
\end{equation}
for all sufficiently small $\tau$.
Suppose the contrary. Then
there exist $\tau_{n}\rightarrow 0$, $x_{n}\in\rit{N}$ $(\|x_{n}\|=1)$ and
$$
y_{n}\in{\absco}[{\class}_{p}(\tau_{n},x_{n})]
$$
such that
$$
y_{n}\rightarrow 0,\qquad ty_{n}\not\in{\absco}[{\class}_{p}(\tau_{n},x_{n})]
\quad {\rm at}\quad t>1\ .
$$
Without loss of generality we can suppose that the sequences $\{x_{n}\}$ and
$\{\frac{y_{n}}{\|y_{n}\|}\}$ are convergent: $x_{n}\rightarrow x$,
$\frac{y_{n}}{\|y_{n}\|}\rightarrow z$.

By Theorem \ref{L11} the linear hull of the set
$\{{\class}_{p}(0,x)\}$ coincides with \rit{N}.
Therefore, there exist matrices
$L_{1}(0)$, $L_{2}(0)$, \ldots , $L_{N}(0)\in{\class}_{p}(0)$ such that
the vectors
$L_{1}(0)x$, $L_{2}(0)x$, \ldots , $L_{N}(0)x$ are linearly independent.
Then the vectors
$L_{1}(\tau_{n})x_{n}$, $L_{2}(\tau_{n})x_{n}$, \ldots , $L_{N}(\tau_{n})x_{n}$
are also linearly independent for all sufficiently large $n$.
Hence, for any positive integer $n$ there exist
$$
\theta ^{(n)}_{1},\theta ^{(n)}_{2}, \ldots , \theta ^{(n)}_{N},\qquad
$$
such that
\begin{equation}\label{au7E}
\sum^{N}_{i=1}\theta ^{(n)}_{i} = 1\ ,
\end{equation}
and the vectors $y_{n}$ are collinear to the respective vectors
\begin{equation}\label{E16}
z_{n} =\sum^{N}_{i=1}\theta ^{(n)}_{i}L_{i}(\tau_{n})x_{n}.
\end{equation}
That is,
\begin{equation}\label{au3E}
z_{n}=\eta _{n}y_{n},
\quad {\rm with}
\quad
\eta _{n}>0\ .
\end{equation}
By definition,  $z_{n}\in{\absco}\{L_{1}(\tau_{n})x_{n}$,
$$
L_{2}(\tau_{n})x_{n}, \ldots , L_{N}(\tau_{n})x_{n}\subseteq
{\absco}\{{\class}_{p}(\tau_{n},x_{n})\}
$$
where $ty_{n}$ does not belong
to the set
${\class}_{p}(\tau_{n},x_{n})$ at $t>1$; hence (\ref{au3E}) implies
$\eta _{n}\le 1$. The last inequality and the condition $y_{n}\rightarrow 0$
imply, in turn,
\begin{equation}\label{au6E}
z_{n}\rightarrow 0.
\end{equation}
The sequences $\{\theta ^{(n)}_{1}\}$, $\{\theta ^{(n)}_{2}\}$, \ldots ,
$\{\theta ^{(n)}_{N}\}$ we can suppose to be convergent to
some limits
$\theta _{1}$, $\theta _{2}$, \ldots , $\theta _{N}$.
As the limit of (\ref{E16}) and (\ref{au7E})
we have:
\begin{equation}\label{au5E}
\sum^{N}_{i=1}\theta _{i}L_{i}(0)x=0,\qquad
{\rm and} \qquad \sum^{N}_{i=1}\theta _{i}=1.
\end{equation}
The relations (\ref{au5E}) contradict the linear independence of the
vectors $L_{1}(0)x$,
$L_{2}(0)x$, \ldots , $L_{N}(0)x$. This contradiction
proves the estimate (\ref{E15}).

\vspace{2mm}
\noindent
Let us return to the proof of the theorem. Denote
$$
\varphi =\liminf_{\tau\rightarrow 0}[{\class}(\tau)],\qquad
\psi =\limsup_{\tau\rightarrow 0}[{\class}(\tau)]
$$
and define  ${\class}_{p}(\tau) =\{L_{1}(\tau)$, $L_{2}(\tau)$, \ldots ,
$L_{Q}(\tau)\}$. Let us establish the inequality
\begin{equation}\label{E17}
\psi\le{\qcm}_{p}[{\class}(0)].
\end{equation}
Chose arbitrary vectors $x\in{\bf S}(1)$, $y\in{\bf S}(\psi )$.
There exist a sequence $\tau_{n}\rightarrow 0$, a sequence
$y_{n}\rightarrow y$ $(y_{n}\in{\bf S}\{{\qcm}_{p}[{\class}(\tau_{n})]\})$
and sequences of real values
$\theta ^{(n)}_{1}$, $\theta ^{(n)}_{2}$, \ldots ,
$\theta ^{(n)}_{Q}$, such that the relations
\begin{equation}\label{E18}
y_{n}=\sum^{N}_{i=1}\theta ^{(n)}_{i}L_{i}(\tau_{n})x,\qquad
\sum^{N}_{i=1}\theta ^{(n)}_{i}\le 1
\end{equation}
hold. Without loss of generality the sequences
$\{\theta ^{(n)}_{1}\}$, $\{\theta ^{(n)}_{2}\}$, \ldots , $\{\theta ^{(n)}_{Q}\}$
can be considered  as convergent:
$$
\theta ^{(n)}_{1}\rightarrow\theta _{1},\quad
\theta ^{(n)}_{2}\rightarrow\theta _{2},\quad \ldots \quad ,\
\theta ^{(n)}_{Q}\rightarrow\theta _{Q}.
$$
Then (\ref{E18}) imply:
\begin{equation}\label{E19}
y=\sum^{N}_{i=1}\theta _{i}L_{i}(0)x,\qquad
\sum^{N}_{i=1}\theta _{i}\le 1.
\end{equation}

Therefore, each vector
$y\in{\bf S}(\psi)$ can be written in the form (\ref{E19})
for any
$x\in{\bf S}(1)$. This proves (\ref{E17}).

Let us establish now the inequality
\begin{equation}\label{E110}
\varphi\ge {\qcm}_{p}[{\class}(0)].
\end{equation}
Because of the inequality
$\varphi\le\psi $, the assertion of the theorem will
follow from  (\ref{E17}) and (\ref{E110}).

If ${\qcm}_{p}[{\class}(0)] = 0$, we have nothing to prove.
Suppose that ${\qcm}_{p}[{\class}(0)] > 0$ and choose some
${\gamma > 0}$ satisfying
\begin{equation}\label{E111}
{\qcm}_{p}[{\class}(0)] - \gamma > 0.
\end{equation}
Let us fix a vector $x\in{\bf S}(1)$ and establish that the condition
\begin{equation}\label{E112}
\|L_{i}(\tau) - L_{i}(0)\|\le
\frac{\gamma\kappa}{{\qcm}_{p}[{\class}(0)]-\gamma }
\end{equation}
with  $\kappa$  from  (\ref{E15}) implies
\begin{equation}\label{E113}
{\bf S}({\qcm}_{p}[{\class}(0)]-\gamma )\subseteq{\absco}{\class}_{p}(\tau)x.
\end{equation}

Let $y$ be an arbitrary vector from
${\bf S}({\qcm}_{p}[{\class}(0)]-\gamma )$. There exist $\theta _{1}$,
$\theta _{2}$, \ldots , $\theta _{Q}$, such that
$$
y=\sum^{Q}_{i=1}
\frac{{\qcm}_{p}[{\class}(0)]-\gamma}
{{\qcm}_{p}[{\class}(0)]}\theta _{i}L_{i}(0)x
$$
and
\begin{equation}\label{E115}
\sum^{Q}_{i=1}\theta _{i}\le 1.
\end{equation}
Hence,
\begin{equation}\label{E116}
y=\sum^{Q}_{i=1}
\frac{{\qcm}_{p}[{\class}(0)]-\gamma}
{{\qcm}_{p}[{\class}(0)]}\theta _{i}L_{i}(\tau)x + z
\end{equation}
where
$$
z=\sum^{Q}_{i=1}
\frac{{\qcm}_{p}[{\class}(0)]-\gamma}
{{\qcm}_{p}[{\class}(0)]}\theta _{i}(L_{i}(0)-L_{i}(\tau))x.
$$
By (\ref{E112}) and (\ref{E115}) the vector $z$ satisfies the estimate
$$
\|z\|\le \frac{\gamma\kappa}{{\qcm}_{p}[{\class}(0)]-\gamma }.
$$
Farther, by (\ref{E15}) there exist $\eta _{1}(\tau)$,
$\eta _{2}(\tau), \ldots , \eta _{Q}(\tau)$
satisfying
\begin{equation}\label{E117}
z=\sum^{Q}_{i=1}
\frac{\gamma}{{\qcm}_{p}[{\class}(0)]-\gamma }\eta _{i}(\tau)L_{i}(\tau)x,
\qquad \sum^{Q}_{i=1}\eta _{i}(\tau)\le 1.
\end{equation}
Define now
\begin{equation}\label{E118}
\theta _{i}(\tau)=
\frac{{\qcm}_{p}[{\class}(0)]-\gamma}
{{\qcm}_{p}[{\class}(0)]}\theta _{i}+
\frac{\gamma}{{\qcm}_{p}[{\class}(0)]-\gamma }\eta _{i}(\tau)\ .
\end{equation}
The relations (\ref{E116}) and (\ref{E117}) imply
$$
y=\sum^{Q}_{i=1}\theta _{i}(\tau)L_{i}(\tau)x\ ,
$$
where
$$
\sum^{N}_{i=1}\theta _{i}(\tau)\le 1
$$
by virtue of
(\ref{E115}), (\ref{E117}),  (\ref{E118}).
We have just proven that the inclusion
(\ref{E113}) holds for all
$\tau$ satisfying (\ref{E112}).
Hence, for such  $\tau$
$$
{\qcm}_{p}[{\class}(\tau)]\ge{\bf S}({\qcm}_{p}[{\class}(0)]-\gamma )\ .
$$
Taking the lower limit of the last inequality at $\tau\rightarrow 0$,
we obtain
$$
\theta\ge{\qcm}_{p}[{\class}(0)]-\gamma .
$$
This and the arbitrariness of  $\gamma > 0$ imply
(\ref{E111}).

The inequalities (\ref{E17}) and  (\ref{E110}) and, consequently, the
theorem are proven.

\section*{Acknowledgements} The authors would like to thank Phil Diamond
for reading a first draft of the manu\-script and making many useful suggestions.

\end{document}